\documentclass[11pt]{amsart}
\usepackage{amssymb, amscd}
\usepackage{eucal}
\usepackage{color}
\usepackage{enumerate}
\usepackage{soul}

\numberwithin{equation}{section}

\DeclareSymbolFont{SY}{U}{psy}{m}{n}
\DeclareMathSymbol{\emptyset}{\mathord}{SY}{'306}

\theoremstyle{plain}

\newtheorem{thm}{Theorem}[section]

\newtheorem{lem}[thm]{Lemma}

\theoremstyle{definition}

\newtheorem{qn}[thm]{Question}
\newcommand{\norm}[1]{\left\lVert#1\right\rVert}
\newcommand{\overbar}[1]{\mkern 1.5mu\overline{\mkern-1.5mu#1\mkern-1.5mu}\mkern 1.5mu}

\newcounter{defcounter}
\title[curvature inequalities]{Curvature inequalities for operators in the \\Cowen-Douglas class of a planar domain}
\author[Md. Ramiz Reza]{Md. Ramiz Reza} 
\address{Department of Mathematics\\Indian Institutte of Science\\Bangalore 560012}
\email{ramiz@math.iisc.ernet.in}

\thanks{This work is supported by the Council of Scientific and Industrial Research (CSIR), India.}
\begin{document}

%\begin{center}
%{\sf\bf Curvature Inequality in Multiply Connected domain}
%\end{center}
%\hrule \vspace*{0.2in}
%
%%---------------------Introduction----------------------
\begin{abstract}
Fix a bounded planar domain $\Omega.$ If an operator $T,$ in the Cowen-Douglas class $B_1(\Omega),$ admits the compact set $\bar{\Omega}$ as a spectral set, then the curvature inequality $\mathcal K_T(w) \leq - 4 \pi^2 S_\Omega(w,w)^2,$ where $S_\Omega$ is the S\"{z}ego kernel of the domain $\Omega,$ is evident. Except when $\Omega$ is simply connected, the existence of an operator for which $\mathcal K_T(w) = 4 \pi^2 S_\Omega(w,w)^2$ for all $w$ in $\Omega$ is not known. However, one knows that if $w$ is a fixed but arbitrary point in $\Omega,$  then there exists a bundle shift of rank $1,$ say $S,$ depending on this $w,$ such that $\mathcal K_{S^*}(w) = 4 \pi^2 S_\Omega(w,w)^2.$  We prove that these {\em extremal}  operators  are uniquely determined: If $T_1$ and $T_2$ are two operators in $B_1(\Omega)$ each of which is the adjoint of a rank $1$ bundle shift and $\mathcal{K}_{T_1}({w}) = -4\pi ^2 S(w,w)^2 = \mathcal{K}_{T_2}(w)$ for a fixed $w$ in $\Omega,$ then $T_1$ and $T_2$ are unitarily equivalent. A surprising consequence is that the adjoint of only some of the bundle shifts of rank $1$ occur as extremal operators in domains of connectivity $> 1.$ These are described explicitly.   
\end{abstract}

\maketitle
\section{Introduction}
Let $\Omega$ be a bounded, open and connected subset of the complex plane $\mathbb C.$  Assume that $\partial \Omega,$ the boundary of $\Omega,$ consists of $n+1$  analytic Jordan curves. Let $\partial \Omega_1, \partial \Omega_2,\cdots, \partial \Omega_{n+1}$ denote the boundary components of $\Omega.$  We shall always let $\partial \Omega_{n+1}$ denote the curve whose interior contains $\Omega.$ Set $\Omega ^{*} =\{\bar{z}\mid z\in \Omega\},$ which is again a planar domain whose boundary consists of $n+1$ analytic Jordan curves. In this paper we study 
operators in $B_1(\Omega^*),$ first introduced by Cowen and Douglas in \cite{CD}, namely, those bounded linear operators $T$ acting on a complex separable Hilbert space $\mathcal H,$ for which $\Omega^* \subseteq \sigma (T)$ and which meet the following requirements.
\begin{enumerate}
   \item ${\rm ran}(T-w) = \mathcal{H},\,\,\,\,  w\in \Omega^*,$
   \item $\bigvee_{w\in \Omega^*} \ker(T-w) = \mathcal{H}$ and
   \item $\dim(\ker(T-w)) = 1,\,\,\,\,  w \in \Omega^*.$ 
\end{enumerate}
These conditions ensure that  one may choose an eigenvector $\gamma_T(w)$ with eigenvalue $w,$ for any operator $T$ in $B_1(\Omega^*),$ such that $w\to \gamma_T(w)$ is holomorphic on $\Omega^*$ (cf. \cite[Proposition 1.11]{CD}). This is the holomorphic frame for the operator $T.$
Cowen and Douglas also provide a model for the operators in the class $B_1(\Omega^*),$ which is easy to describe: 

If $T\in B_1(\Omega^*)$ then $T$ is unitarily equivalent to the adjoint  $M^*$ of the operator of multiplication $M$ by the coordinate function on some Hilbert space $\mathcal{H}_{{K}}$ consisting of holomorphic function on $\Omega$ possessing a reproducing kernel $K$. Throughout this paper, we let $M$ denote the operator of multiplication by the coordinate function and as usual  $M^*$ denotes its adjoint.

The kernel $K$ is complex valued function defined on $\Omega\times \Omega,$ which is holomorphic in the first and anti-holomorphic in the second variable and is positive definite in the sense that 
$
\big (\!\!\big ( K(z_i,z_j) \big )\!\! \big )
$
is positive definite for every subset $\{z_1,\ldots,z_n\}$ of the domain $\Omega.$
We will therefore assume, without loss of generality, that an operator $T$ in $B_1(\Omega^*)$ has been realized as the operator $M^*$ on some reproducing kernel Hilbert space $\mathcal H_K.$ The curvature $\mathcal K_T$ of the operator $T$ is defined as 
%\begin{align*}
$$ \mathcal{K}_T(\bar{z}) = -\frac{\partial ^2}{\partial w \partial \bar{w}}\mbox{log} K_{T}(w,w)|_{w=z} = -\frac{\norm{K_{z}}^2 \norm{\bar{\partial}K_{z}}^2 - |{\langle K_{z}, \bar{\partial}K_{z} \rangle}|^2}{(K({z},{z}))^2},\;\;z \in \Omega, $$ 
%\end{align*}
where $K_z$ and $\bar{\partial}K_z$  are the vectors 
\begin{align*}
 K_z(u) &:= K(u,z)\;\;,u\in \Omega, \\ \bar{\partial}K_z (u) &:= \frac{\partial}{\partial \bar{w}}K(u,w)\mid _{w=z}\;\;,u\in \Omega,
\end{align*}
in $\mathcal H_K.$ Thus the curvature is a real analytic function on $\Omega^*.$ 

It turns out that this definition of curvature is independent of the representation of the operator $T$ as the adjoint $M^*$ of a  multiplication operator $M$ on some reproducing kernel Hilbert space 
$\mathcal{H}_K.$  Indeed, if $T$ also admits a representation as the adjoint of the multiplication operator on another reproducing kernel Hilbert space $\mathcal H_{\tilde{K}},$ then we must have $K(z,w) = \varphi(z)\tilde{K}(z,w) \overline{\varphi(w)}$ for some holomorphic function $\varphi$ defined on $\Omega.$(See \cite[Section 1.15]{CD}). This implies $\frac{\partial ^2}{\partial w \partial \bar{w}}\mbox{log} K_{T}(w,w)|_{w={z}}  = \frac{\partial ^2}{\partial w \partial \bar{w}}\mbox{log} \tilde{K}_{T}(w,w)|_{w={z}}.$ 

If $T_1$ and $T_2$ are any two operators in $B_1(\Omega),$ then any intertwining unitary must map the holomorphic frame of one to the other.  From this, it follows that the curvatures of these two operators must be equal, as shown in \cite[Theorem 1.17]{CD} along with the non-trivial converse.	
\begin{thm}
Two operators $T_1$ and $T_2$ in $B_1(\Omega ^*)$ are unitarily equivalent if and only if their associated curvature functions are equal that is $\mathcal{K}_{T_1}(w) = \mathcal{K}_{T_2}(w)$ for all $w \in \Omega^*.$ 
\end{thm}

%We will assume through out this paper that the operator $T$ is in $B_1(\Omega ^*).$ and $K_T(z,w)$ be a associated Reproducing kernel on $\Omega \times \Omega$, so that $T$ is unitarily equivalent to  $M_z^*$, the adjoint of Multiplication operator on the Reproducing kernel Hilbert space $\mathcal{H}(K_T)$. Note that as $K_T(z,w)$ is holomorphic in first variable and anti holomorphic in second variable, it turns out that Reproducing kernel Hilbert space $\mathcal{H}(K_T)$ consist of holomorphic function on $\Omega$. Now the curvature of the operator $T$ is defined as follows
%\begin{align}\label{eq:Curvature definition}
%\mathcal{K}_T(z) = -\frac{\partial ^2}{\partial w \partial \bar{w}}\mbox{log} K_{T}(w,w)|_{w=\bar{z}}, \;\;\mbox{for all}\;\;z \in \Omega ^*
%\end{align}

Recall that a compact subset $X\subseteq \mathbb C$ is said to be a spectral set for an operator $A$ in $\mathcal{L}(\mathcal{H}),$  if

\begin{equation*}
\sigma (A) \subseteq X\;\;\mbox{and}\;\;\sup \{||r(A)|| \mid r\in \mbox{\rm{Rat}}(X)\;\;\mbox{and}\;\; ||r||_{\infty} \leq 1\} \leq 1, 
\end{equation*} 
where $\mbox{ \rm Rat}(X)$ denotes the algebra of rational function whose poles are off $X$ and $||r||_{\infty} $ denotes the sup norm over the compact subset $X.$  Equivalently, $X$ is a spectral set for the operator $A$ if the homomorphism $\rho_A: \mbox{\rm Rat}(X) \to \mathcal L(\mathcal H)$ defined by the formula  $\rho_A(r) = r(A)$ is contractive. There are plenty of examples where the spectrum of an operator is a spectral set, for instance, this is the case for  subnormal operators (See \cite[Chapter 21]{HALbook}).

Now assume  $\overbar{\Omega^*},$ the closure of $\Omega ^*,$ is a spectral set for the operator $T$ in $B_1(\Omega^*).$ The space $\ker(T-w)^2 = \rm span\{K_w, \bar{\partial}K_w\}$ is an invariant subspace for $T.$ Representing the restriction of the operator $T$ to this subspace with respect to an orthonormal basis as a $2\times 2$ matrix, we have 
$$
T_{|\ker(T-w)^2} = \Big ( \begin{smallmatrix} w & \tfrac{1}{\sqrt{- \mathcal K_T(w)}}\\ 0 & w\end{smallmatrix} \Big ).  
$$
It follows that $\overbar{\Omega ^*}$ is also a spectral set for $T\mid _{\rm ker(T-w)^2}.$ For any $r$ in $\mbox{\rm Rat}(\overline{\Omega^*}),$ it is not hard to verify that 
$$r\big (T_{|\ker(T-w)^2}\big ) = \Big ( \begin{smallmatrix} r(w) &  \tfrac{r^\prime(w)}{\sqrt{- \mathcal K_T(w)}}\\ 0 & r(w)\end{smallmatrix} \Big ).  
$$
Since 
\begin{align*}
\sup \{|r^\prime(w)|\mid \|r\|_\infty \leq 1,\,r\in \mbox{\rm Rat}(\overline{\Omega^*})\} = 
2\pi (S_{\Omega ^*}(w,w)),\,w\in \Omega^*,
\end{align*}
where $S_{\Omega ^*}(z,w),$ the S\"{z}ego kernel of $\Omega^*,$ is the reproducing kernel for the Hardy space $(H^2(\Omega ^*),ds),$ a curvature inequality becomes evident (see \cite[Corollary 1.2]{GMCI}), that is, 
\begin{align}\label{eq:Curvature Inequality}
\mathcal{K}_{T} (w) &\leq -4\pi ^2(S_{\Omega ^*}(w,w))^2,\;\;\;w\in\Omega ^*.
\end{align}
% This inequality we call as curvature inequality for an operator in $B_1(\Omega ^*)$ which possessed $\overline{\Omega ^*}$ as a spectral set.
Equivalently, since $S_{\Omega }(z,w) = S_{\Omega^* }(\bar{w},\bar{z}),$ the curvature inequality takes the form 
\begin{align}\label{eq: alternative CI}
\frac{\partial ^2}{\partial w \partial \bar{w}}\mbox{log} K_{T}(w,w)\geq 4\pi ^2 (S_{\Omega }(w,w))^2, \;\;\;w \in \Omega.
\end{align}

%It is known that Szego kernel for the domain $\Omega$ and the associated %conjugate domain $\Omega ^*$ are related by following relation
%\begin{align}\label{eq: relation between Szegos}
%S_{\Omega }(z,w) &= S_{\Omega^* }(\bar{w},\bar{z})
%\end{align}
%Hence using $\eqref{eq: relation between Szegos}$, the curvature inequality can also be stated as follows

The operator $M^*$ on the  Hardy space $(H^2(\mathbb{D}), ds)$ is in $B_1(\mathbb{D} ).$  The closed unit disc  is a spectral set for the operator $M^*.$ The reproducing kernel of the Hardy space, as is well-known, is the  S\"{z}ego kernel $S_{\mathbb{D}}$ of the unit disc $\mathbb D.$ It is given by the formula  $S_{\mathbb{D}}(z,w) = \frac{1}{2\pi(1-z\bar{w})},$ for all $z,w$ in $\mathbb{D}.$  The computation of the curvature of the operator $M^*$ is  straightforward and is given by the formula 
\begin{align*}
-\mathcal{K}_{M^*} (w)= \frac{\partial ^2}{\partial w \partial \bar{w}}\mbox{log} S_{\mathbb{D} }(w,w) &= 4 \pi ^2 (S_{\mathbb{D}}(w,w))^2,\;\;\;w\in \mathbb D.
\end{align*}
Since the closed unit disc is a spectral set for any contraction, it follows that the curvature of the operator $M^*$ on the Hardy space $(H^2(\mathbb{D}), ds)$  dominates the curvature of every other contraction in $B_1(\mathbb D).$ 

If the region $\Omega$ is simply connected, then using the Riemann map and the transformation rule for the S\"{z}ego kernel (See \cite[Theorem 12.3]{BELLcauchy}) together with  the chain rule for composition,  we see that 
\begin{align}\label{eq: Szego in simply connected}
\frac{\partial ^2}{\partial w \partial \bar{w}}\mbox{log} S_{\Omega }(w,w) &= 4 \pi ^2 (S_{\Omega}(w,w))^2,\;\;w\in \Omega.
\end{align}
% This tells us that  the adjoint $M_z^*$ of multiplication operator on $(H^2(\Omega),ds)$ the usual Hardy space on $\Omega $. It is straightforward to verify that $M_z ^*$ is in $B_1(\Omega ^*)$ and $\overline{\Omega ^*}$ is a spectral set for $M_z^*$. As a consequence of $\eqref{eq: Szego in simply connected}$, It turns out that
% \begin{align*}
% \mathcal{K}_{M_z ^*} (w) &= -4\pi ^2(S_{\Omega ^*}(w,w))^2,\;\;\mbox{for all}\;w %\in \Omega ^ *.
% \end{align*}
This shows that in the case of bounded simply connected domain with jordan analytic boundary, the operator $M^*$ on $(H^2(\Omega),ds)$ is an extremal operator.
 
On the other hand, if the region is not simply connected, then \eqref{eq: Szego in simply connected} fails. Indeed, Suita (see \cite{SUITA}) has shown that 
\begin{align} \label{eq: not extremal}
\frac{\partial ^2}{\partial w \partial \bar{w}}\mbox{log} S_{\Omega }(w,w) &> 4 \pi ^2 (S_{\Omega}(w,w))^2,\;\;w\in \Omega. 
 \end{align}
 Or equivalently,
 \begin{align}\label{eq: not extremal2}
 \mathcal{K}_{M^*} (z) &< -4\pi ^2(S_{\Omega ^*}(z,z))^2,\;\;z \in \Omega ^*,
 \end{align}
where $M^*$ is the adjoint of the multiplication by the coordinate function on the Hardy space $(H^2(\Omega),ds).$ We therefore conclude that if $\Omega$  is not simply connected, then the operator $M^*$ fails to be extremal.
 
We don't know if there exists an operator $T$ in $B_1(\Omega ^*),$ admitting  
$\overbar{\Omega ^*}$ as a spectral set for which 
 \begin{align*}
\mathcal{K}_{T} (w) = -4\pi ^2(S_{\Omega ^*}(w,w))^2,\;\;w \in \Omega ^*.
 \end{align*}
The question of equality at just one fixed but arbitrary point $\bar{\zeta}$ in $\Omega^*$
was answered in \cite[Theorem 2.1]{GMCI}. An operator  $T$ in $B_1(\Omega ^*),$ which admits  $\overline{\Omega^*}$ as a spectral set would be called extremal at $\bar{\zeta}$ if  
\begin{align*}
 \mathcal{K}_{T} (\bar{\zeta}) &= -4\pi ^2(S_{\Omega ^*}(\bar{\zeta},\bar{\zeta}))^2.
\end{align*} 
Equivalently, representing the operator $T$ as the operator $M^*$ on a Hilbert space  possessing a reproducing kernel $K_T: \Omega\times \Omega \to \mathbb C,$ we have that 
\begin{equation} \label{eq:equality at zeta}
 \frac{\partial ^2}{\partial w \partial \bar{w}}\mbox{log} K_{T}(w,w)|_{w=\zeta} = 4 \pi ^2 S_{\Omega}(\zeta,\zeta)^2. 
\end{equation}

Since the the operator $M^*$ on the Hardy space $H^2(\mathbb D)$ is extremal, that is, $ \mathcal K_T(w) \leq \mathcal K_{M^*}(w),$ for all $w\in \mathbb D$ and for every $T$ in the class of contractions in $B_1(\mathbb D)$ and the curvature is a complete invariant, one may ask: 
\begin{qn}[R. G. Douglas]
If $\mathcal K_T(w_0) = \mathcal K_{M^*}(w_0)$ for some fixed $w_0$ in $\mathbb D,$ then does it follow that $T$ must be unitarily equivalent to $M^*$?
\end{qn} 
This question has an affirmative answer if, for instance, $T$ is a homogeneous operator and it is easy to construct examples where the answer is negative (cf. \cite{GMback}). 

The operator $M$ on the Hardy space $H^2(\mathbb D)$  is a pure subnormal operator with the property: the spectrum of the minimal normal extension,  designated  the normal  spectrum, is contained in the boundary of the spectrum of the operator $M.$ These properties determine the operator $M$ uniquely up to unitary equivalence. The question of characterizing all pure subnormal operator with spectrum $\overline{\Omega}$ and normal spectrum contained  in the boundary of $\Omega$ is more challenging if $\Omega$ is not simply connected.  The deep results of Abrahamse and Douglas (cf. \cite[Theorem 11]{ABsubnormal})  show that these are exactly the bundle shifts, what is more,  they are in one to one correspondence with the equivalence classes of flat unitary bundles on the domain $\Omega.$ It follows that adjoint of a bundle shift of rank $1$ lies in $B_1(\Omega ^*).$ Since bundle shifts are subnormal with spectrum equal to $\overline{\Omega}, $ it follows that $\overline{\Omega ^*}$  is a spectral set for the adjoint of the bundle shift. In fact, the extremal operator at $\bar{\zeta},$  found in \cite{GMCI}, is the adjoint of a bundle shift of rank $1$. Therefore, one may ask, following R. G. Douglas, if the curvature $\mathcal K_T(\bar{\zeta})$  of an operator $T$ in $B_1(\Omega^*),$ admitting $\overline{\Omega^*}$ as  a spectral set, equals $-4 \pi^2 S_{\Omega}(\zeta,\zeta)^2,$ then does it follow that $T$ is necessarily unitarily equivalent to the extremal operator at $\bar{\zeta}$ found in \cite{GMCI}. In this paper, we show that an extremal operator must be uniquely determined  within  $\big \{[\![T^*]\!]: T \mbox{\rm\,\,  is bundle shift of multiplicity } 1 \mbox{ over }\Omega \big \},$ where $[\![\cdot]\!]$ denotes  the unitary equivalence class.

\section{Preliminaries on bundle shifts of rank one}
%In this section we will recall some basic well known facts about Bundle shift. 
Let $\alpha$ be an element in $\rm Hom(\pi_1(\Omega), \mathbb{T})$ that is it is a homomorphism from the fundamental group $\pi_1(\Omega)$ of $\Omega$ into the unit circle $\mathbb{T}.$ Such homomorphism is also called a character. Each of these character induce a flat unitary bundle $E_\alpha$ of rank $1$ on $\Omega$ (cf. \cite[Proposition 2.5]{CHENdual}). Following theorem establishes one to one correspondence between $\rm Hom(\pi_1(\Omega), \mathbb{T})$ and the set equivalence classes of flat unitary vector bundle over $\Omega$ of rank $1$(See \cite[p. 186]{GUNNINGrs}.)

%Each such $\alpha$ induce a rank one flat unitary vector bundle $E_{\alpha}$ over the domain $\Omega.$ In fact every rank one flat unitary vector bundle over the domain $\Omega$ is isomorphic to one such flat unitary vector bundle $E_{\alpha}$ for some homomorphism $\alpha.$ The following well known theorem characterize the equivalence class of rank one flat unitary vector bundle.
\begin{thm}
Two rank one flat unitary vector bundle $E_{\alpha}$ and $E_{\beta}$ are equivalent as flat unitary  vector bundle if and only if their inducing characters are equal that is $\alpha = \beta.$
\end{thm}

First, if $f$ is a holomorphic section of the bundle $E_\alpha,$ then for $z\in U_i\cap U_j,$ where $\{U_i, \phi_i\}_{i\in I}$ is a coverng of $\Omega,$ we have that $|(\phi _i ^z)^{-1} (f(z))| = |(\phi _j^ z)^{-1} (f(z))|.$
Thus the function $h_f(z) := |(\phi_i^z)^{-1}(f(z))|,$ $z\in U_i,$ is well defined on all of $\Omega$  and is subharmonic there.  
Let $H^2_{E_{\alpha}}$ be the linear space of those holomorphic sections 
$f$ of $E_{\alpha}$ such that the subharmonic function $(h_f)^2$  on $\Omega$ is majorized by a harmonic function on $\Omega.$
While there is no natural inner product on  the space $H^2_{E_{\alpha}},$  
Abrahamse and Douglas define an inner product relative to the harmonic measure with respect to a fixed but arbitrary point $p\in \Omega.$ We make the comment in \cite[p. 118]{ABsubnormal} explicit in what follows. 
Let $\{\Omega_k\}_{k \in \mathbb{N}}$ be a regular exhaustion $\Omega$ that is,  it is a sequence of increasing subdomains $\Omega_k$ of $\Omega$  satisfying
\begin{enumerate}[(a)]
\item $\overbar{\Omega}_k \subset \Omega_{k+1},$
\item $\cup_k \Omega_k = \Omega,$
\item $p \in \Omega_1$
\item boundary of each $\Omega_k$ consists of finitely many smooth jordan curve.
\end{enumerate}
Then norm of the section $f$ in $H^2_{\alpha}(\Omega)$ is then defined by the limit 
\begin{align*}
||f||^2 = \lim _{k \to \infty}-\frac{1}{2\pi}\int _{\partial \Omega_k} (h_f(z))^2\frac{\partial}{\partial \eta _z}\big(g_{k}(z,p)\big)ds(z),
\end{align*} where $g_{k}(z,\zeta)$ denote the green function for the domain $\Omega_k$ at the point $p$ and $\frac{\partial}{\partial \eta _z} $ denote the directional derivative along the outward normal direction w.r.t the positively oriented boundary of $\Omega_k.$ The linear space $H^2_{E_{\alpha}}$ is complete with respect to this norm making it into a Hilbert space.  A bundle shift $T_{E_{\alpha}}$ is simply the operator of multiplication  by the coordinate function on  $H^2_{E_{\alpha}}.$ 
\begin{thm} [Abrahamse and Douglas] 
Let $E_{\alpha}$ and $E_{\beta}$ be rank one flat unitary vector bundles induced by the  homomorphisms $\alpha$ and $\beta$ respectively. Then the bundle shift $T_{E_{\alpha}}$ is unitarily equivalent to the bundle shift $T_{E_{\beta}}$ if and only if  $E_{\alpha}$ and $E_{\beta}$ are equivalent as flat unitary vector bundles.  
\end{thm}

It is not very hard to verify that $T_{E_{\alpha}}$ is a pure cyclic subnormal operator with spectrum  $\overbar{\Omega}$ and normal spectrum  $\partial \Omega.$ 
In fact these are the characterizing property for the rank one bundle shift.
%Abrahamse and Douglas proved that above property are the characterizing property for Multiplicity one Bundle shift, namely
\begin{thm}[Abrahamse and Douglas]
Every pure cyclic subnormal operator with spectrum  $\overline{\Omega}$ and normal spectrum contained in $\partial \Omega$ is unitarily equivalent to a bundle shift $T_{E_{\alpha}}$ for some character $\alpha.$
\end{thm}

Bundle shifts can also be realized as a multiplication operator on a certain subspace of the classical Hardy space $H^2(\mathbb{D}).$ Let $\pi:\mathbb{D}\mapsto \Omega$ be a holomorphic covering map satisfying $\pi (0) = p.$ Let $G$ denote the group of deck transformation associated to the map $\pi$ that is $G = \{A \in \rm Aut(\mathbb{D}) \mid \pi \circ A = \pi \}.$ As $G$ is isomorphic to the fundamental group $\pi _1(\Omega)$ of $\Omega$, every character $\alpha$ induce a unique element in $\rm Hom(G, \mathbb{T}).$ By an abuse of notation we will also denote it by $\alpha.$ A  holomorphic function $f$ on unit disc $\mathbb{D}$ satisfying $f \circ A = \alpha (A) f , \;\;\text{for all}\; A \in G$, is called a modulus automorphic function of index $\alpha.$ Now consider the following subspace of the Hardy space $H^2(\mathbb{D})$ which consists of modulus automorphic function of index $\alpha$, namely  
\begin{align*}
H^2(\mathbb{D},\alpha) = \{ f \in H^2(\mathbb{D}) \mid f \circ A = \alpha (A) f , \;\;\text{for all}\; A \in G\}
\end{align*}
 Let $T_{\alpha}$ be the multiplication operator by the covering map $\pi$ on the subspace $H^2(\mathbb{D},\alpha).$ Abrahamse and Douglas have shown in \cite[Theorem 5]{ABsubnormal} that the operator $T_{\alpha}$ is unitarily equivalent to the bundle shift $T_{E_{\alpha}}.$  
 
There is another way to realize the bundle shift as a multiplication operator $M$ on a Hilbert space of multivalued holomorphic function defined on  $\Omega$ with the property that its absolute value is single valued. A multivalued holomorphic  function defined on  $\Omega$ with the property that its absolute value is single valued is called a multiplicative function. Every modulus automorphic function $f$ on $\mathbb{D}$ induce a multiplicative function on $\Omega,$ namely, $f \circ \pi ^{-1}.$ Converse is also true (See \cite [Lemma 3.6]{VOICHICKideal}). We define the class $H^2_{\alpha}(\Omega)$  consisting of multiplicative function in the following way: 
\begin{align*}
H^2_{\alpha}(\Omega):=\{f\circ \pi^{-1}\mid f\in H^2(\mathbb{D},\alpha) \}
\end{align*}

So the linear space $H^2_{\alpha}(\Omega)$ is consisting of those multiple valued function $h$ on $\Omega$ for which $|h|$ is single valued, $|h|^2$ has a harmonic majorant on $\Omega$ and $h$ is locally holomorphic in the sense that each point $w \in \Omega$ has a neighbourhood $U_w$ and a single valued holomorphic function $g_{w}$ on $U_w$ with the property $|g_{w}| = |h|$ on $U_w$ (See \cite[p.101]{FISHER}). Since the covering map $\pi$ lifts the harmonic measure $d\omega_p$ on $\partial \Omega$ at the point $\pi (0)= p$  to the linear Lebesgue measure on the unit circle $\mathbb{T}$, It follows that  $H^2_{\alpha}(\Omega)$ endowed with the norm
\begin{align*}
||f||^2 &= \int _{\partial \Omega} |f(z)|^2 d\omega_p(z),                        
\end{align*} becomes a Hilbert space (cf. \cite[p. 101]{FISHER}.) We will denote it by $\big(H^2_{\alpha}(\Omega),d\omega_{p}\big).$ In fact the map $f \mapsto f \circ \pi ^{-1}$ is a unitary map from $H^2(\mathbb{D},\alpha)$ onto $\big(H^2_{\alpha}(\Omega),d\omega_{p}\big)$ which intertwine the multiplication by $\pi$ on $H^2(\mathbb{D},\alpha)$ and the multiplication by coordinate function $M$ on $\big(H^2_{\alpha}(\Omega),d\omega_{p}\big).$

We have described three different but unitarily equivalent realization of a bundle shift of rank $1$ over the domain $\Omega.$ We prefer to work with the third realization. It is well known that the harmonic measure $d\omega_{p}$ on $\partial \Omega$  at the point $p$ is boundedly mutually absolutely continuous w.r.t the arc length measure $ds$ on $\partial \Omega.$ In fact we have 
\begin{align*}
d\omega _{p}(z) &= -\frac{1}{2\pi} \frac{\partial}{\partial \eta _z}\big(g(z,p)\big)ds(z),\;\; z \in \partial \Omega, 
\end{align*} where $g(z,\zeta)$ denote the green function for the domain $\Omega$ at the point $p$ and $\frac{\partial}{\partial \eta _z} $ denote the directional derivative along the outward normal direction (w.r.t positively oriented $\partial \Omega$). In this paper, instead of working with harmonic measure $d\omega_{p}$ on $\partial \Omega$,  we will work with arclength measure $ds$ on $\partial \Omega.$ This is the approach in Sarason \cite{SARA}. So, we define the norm of a function $f$ in in $H^2_{\alpha}(\Omega)$ by
\begin{align*}
||f||^2_{ds} = \int _{\partial \Omega} |f(z)|^2ds 
\end{align*}
Since the outward normal derivative of the Green's function is negative on the boundary $\partial \Omega,$ we have 
\begin{align}\label{eq: harmonic}
d\omega _{p}(z) = h^2(z) ds(z),\;\;z \in \partial \Omega,
\end{align}where $h(z)$ is a positive continuous function on $\partial \Omega.$ We also see that 
\begin{align*}
c_2||f||^2_{ds} \leq ||f||^2 \leq c_1||f||^2_{ds},
\end{align*}
where $c_1$ and $c_2$ are the supremum and the infimum of the function $h$ on $\partial \Omega.$ 

Hence it is clear that $||\cdot||_{ds}$ defines an equivalent norm on   $H^2_{\alpha}(\Omega),$ We let $\big(H^2_{\alpha}(\Omega),ds\big)$ be the Hilbert space which is the same as $H^2_{\alpha}(\Omega)$ as a linear space but is given the new norm $||\cdot||_{ds}.$  In fact, the identity map from $\big(H^2_{\alpha}(\Omega),d\omega_{p}\big)$ onto $\big(H^2_{\alpha}(\Omega),ds\big)$ is invertible and intertwines the corresponding multiplication operator by the coordinate function. It is easily verified that the multiplication operator by coordinate function on $\big(H^2_{\alpha}(\Omega),ds\big)$ is also a pure  cyclic subnormal operator with spectrum equal to $\overline{\Omega}$ and normal spectrum contained  in $\partial \Omega.$ By a slight  abuse of notation, we will denote the multiplication operator by the coordinate function on $\big(H^2_{\alpha}(\Omega),ds\big)$ also by $T_{\alpha}.$  

Using the characterization of all cyclic subnormal operator with spectrum equal to $\overline{\Omega}$ and normal spectrum contained in $\partial \Omega$ given by Abrhamse and Douglas, we conclude that for every character $\beta,$ the operator $T_{\beta}$ on  $\big(H^2_{\beta}(\Omega),ds\big)$ is unitarily equivalent to $T_{\alpha}$ on $\big(H^2_{\alpha}(\Omega),d\omega_p\big)$ for some $\alpha.$ 
In the following section we will establish a bijective correspondence  (which respects the unitary equivalence class) between these two kinds of bundle shifts. The following Lemma helps in establishing this bijection. 

\begin{lem}\label{existence of modulus auto}
If $v$ be a positive continuous function on $\partial \Omega,$ then there exist a  function $F$ in $H^{\infty}_{\gamma}(\Omega)$ for some character $gamma$ such that $|F|^2 = v$ almost everywhere (w.r.t  arc length measure), on $\partial \Omega.$ In fact $F$ is invertible in the sense that there exist $G$ in $H^{\infty}_{\gamma ^{-1}}(\Omega)$ so that $FG = 1$  on $\Omega.$
\end{lem}
\begin{proof}
Since $v$ is a positive continuous function on  $\partial \Omega,$ it follows that $\log v$ is continuous on $\partial \Omega.$ Since the boundary $\partial \Omega$ of $\Omega$ consists of jordan analytic curves, the Dirichlet problem is solvable with continuous boundary data. Now solving the Dirichlet problem with boundary value $\frac{1}{2}\log v,$ we get a harmonic function $u$ on $\Omega$ with continuous boundary value $\frac{1}{2}\log v.$ Let $u^*$ be the multiple value conjugate harmonic function of $u.$ Let's denote the period of the multiple valued conjugate harmonic function $u^*$ around the boundary component $\partial \Omega_{j}$ by 
\begin{align*}
   c_j = -\int _{\partial \Omega _{j}} \frac{\partial}{\partial \eta _z}\big(u(z)\big)ds_z, \;\;\;\mbox{for}\; j= 1,2,...,n\;\;
\end{align*}
In the above equation negative sign appear since we have assumed that $\partial \Omega$ is positively oriented, hence the different components of the boundary $\partial \Omega_j,$ $j=1,2,\ldots, n,$ except the outer one are  oriented in clockwise direction. Now consider the function $F(z)$ defined by 
\begin{align*}
F(z) &= \exp(u(z) + i u^*(z)) 
\end{align*} 
Now observe that $F$ is a multiplicative holomorphic  function  on $\Omega.$ Hence following \cite[Lemma 3.6]{VOICHICKideal}, we have a existence of modulus automorphic function $f$ on unit disc $\mathbb{D}$ so that $F = f \circ \pi ^{-1}.$ We find the index of the modulus automorphy for the function $f$ in the following way. Around each boundary component $\partial \Omega_{j},$ along the anticlockwise direction, the value of $F$ gets changed by $\exp(i c_j)$ times its initial value. So, the index of $f$ is determined by the $n$ tuple of numbers $(\gamma_1, \gamma_2,...,\gamma_n)$ given by,
\begin{align*}
\gamma_j &= \exp (i c_j),  \;\;\;\; j= 1,2,...,n.\;\;
\end{align*}
For each of these $n$ tuple of numbers, there exist a homomorphism $\gamma :\pi_1(\Omega) \to \mathbb{T}$ such that these $n$ tuple of numbers occur as  a image of the $n$ generator of the group $\pi_1(\Omega)$ under the map $\gamma.$
Also we have $|F(z)|^2 = \exp (2u(z)) = v(z), \;\;\;z \in \partial \Omega.$ Since  $u$ is continuous on $\overline{\Omega}$, it follows that  $|F(z)|$ is bounded on $\Omega.$ Hence $F$ belongs to $H^{\infty}_{\gamma}(\Omega)$ with $|F|^2 =v$ on $\partial \Omega.$ 

The function $\frac{1}{v}$ is also positive and continuous on $\partial \Omega$, as before, there exists a function $G$ in $H^{\infty}_{\delta}(\Omega)$ with $|G|^2 =\frac{1}{v}$ on $\partial \Omega.$ Since $\log\frac{1}{v} = - \log v$, it easy to verify that index of $G$ is exactly $(\gamma_1^{-1}, \gamma_2^{-1},...,\gamma_n^{-1})$ and hence $\delta$ is equal to $\gamma^{-1}.$ Evidently $FG =1$ on $\Omega.$
\end{proof}

Now we establish the bijective correspondence which preserve the unitary equivalence class, promised earlier. 
From \eqref{eq: harmonic}, we know that the harmonic measure $d\omega_p$ is of the form $h^2ds$ for some positive continuous function on $\partial \Omega.$ Combining this with the preceding Lemma, we see that there is a $F$ in $H^{\infty}_{\gamma}(\Omega)$ with $|F|^2 = h^2$ on $\partial \Omega$ and a $G$ in $H^{\infty}_{\gamma ^{-1}}(\Omega)$ with $|G|^2 = h^{-2}$ on $\partial \Omega.$ Now consider the map  $M_F:\big(H^2_{\alpha}(\Omega),d\omega_p(z)\big)\mapsto \big(H^2_{\alpha \gamma }(\Omega),ds\big),$ defined by the equation 
\begin{align*}
%\;\;\; M_F:\big(H^2_{\alpha}(\Omega),d\omega_p(z)\big)\mapsto \big(H^2_{\alpha \gamma }(\Omega),ds\big),\;\;\mbox{defined by,}\\
M_F(g) =F g, \;\;\;\;g\in \big(H^2_{\alpha}(\Omega),d\mu_{\zeta}(z)\big). 
\end{align*}
Clearly, $M_F$ is a unitary operator and its inverse is the operator $M_{G}.$ The multiplication operator $M_F$ intertwines the corresponding operator of multiplication by the coordinate function on the Hilbert spaces $\big(H^2_{\alpha}(\Omega),d\omega_p(z)\big)$ and $\big(H^2_{\alpha \gamma }(\Omega),ds\big)$
establishing a bijective correspondence of the unitary equivalence classes of bundle shifts.  As a consequence we have the following theorem which was  proved by Abrahamse and Douglas (See \cite[Theorem 5 and 6]{ABsubnormal}) with the harmonic measure $d\omega_p$ instead of the arc length measure $ds.$
 \begin{thm}\label{Equivalence of character}
The bundle shift $T_{\alpha}$ on $\big(H^2_{\alpha}(\Omega),ds\big)$ is unitarily equivalent to the bundle shift $T_{\beta}$ on $\big(H^2_{\beta}(\Omega),ds\big)$ iff $\alpha = \beta.$ 
\end{thm}

It can be shown using the result of Abrahamse and Douglas (See \cite[Theorem 3]{ABsubnormal}) that for any character $\alpha$, the adjoint of the rank $1$ bundle shift $T_{\alpha}$ lies in $B_1(\Omega^*).$ Since the bundle shifts $T_{\alpha}$ is subnormal, it follows that the adjoint of the bundle shifts $T_{\alpha}$ admits $\Omega ^*$ as a spectral set. Consequently, we have an inequality for the curvature of the bundle shifts, namely, 
\begin{align*}
\mathcal{K}_{T_{\alpha}^*} (w) \leq -4\pi ^2(S_{\Omega ^*}(w,w))^2,\;w  \in \Omega ^*.
\end{align*}   
Given any fixed but arbitrary point $\zeta$ in $\Omega,$ in the following section, we  recall the proof (slightly different from the original proof given in \cite{GMCI} of the existence of a bundle shift $T_{\alpha}$ for which equality occurs at $\bar{\zeta}$ in the curvature inequality. However, the main theorem of this paper is the ``uniqueness'' of such an operator. 
\begin{thm}[Uniqueness]\label{Uniqueness}
If the bundle shift $T_{\alpha}$ on $\big(H^2_{\alpha}(\Omega),ds\big)$ and the bundle shift $T_{\beta}$ on $\big(H^2_{\beta}(\Omega),ds\big)$ are extremal at the point $\bar{\zeta},$ that is, if they satisfy 
\begin{align*}
\mathcal{K}_{T_{\alpha}^*} (\bar{\zeta}) &= -4\pi ^2(S_{\Omega ^*}( \bar{\zeta},\bar{\zeta}))^2 = \mathcal{K}_{T_{\beta}^*} (\bar{\zeta})
\end{align*}  then the bundle shifts $T_{\alpha}$ and $T_{\beta}$ are unitarily equivalent, which is the same as  $\alpha = \beta.$
\end{thm}

The  Hardy space $\big( H^2(\Omega), d\omega_p\big)$ consists of holomorphic function on $\Omega$ such that $|f|^2$ has a harmonic majorant on $\Omega.$ Each $f$ in $\big( H^2(\Omega), d\omega_p\big)$ has a non tangential boundary value almost everywhere. In the usual way $\big( H^2(\Omega), d\omega_p\big)$ is identified with a closed subspace of $L^2(\partial \Omega, d\omega_p)$ (see \cite[Theorem 3.2]{RUDINhp}). Let $\lambda$ be a positive continuous function on $\partial \Omega$. As the measure $\lambda ds$ and the harmonic measure $d\omega_p$ on $\partial\Omega$ are boundedly mutually absolutely continuous one can define an equivalent norm on $\big( H^2(\Omega)$ in the following way
\begin{align*}
||f||^2_{\lambda ds} &= \int _{\partial \Omega} |f(z)|^2 \lambda(z) ds(z).
\end{align*}
Let $\big( H^2(\Omega), \lambda ds\big)$ denote the linear space $ H^2(\Omega)$ endowed with the  norm $\lambda ds.$ Since the harmonic measure $d\omega_p$ is boundedly mutually absolutely continuous with the arc length measure and $\lambda$ is a positive continuous function on $\partial \Omega,$ it follows that identity map $id: \big( H^2(\Omega), d\omega_p\big) \mapsto \big( H^2(\Omega), \lambda ds\big)$ is an invertible map intertwining the associated multiplication operator $M.$ Thus $\big( H^2(\Omega), \lambda ds\big)$ acquires the structure of a Hilbert space and the operator $M$ on it is cyclic,  pure subnormal, its spectrum is equal to $\overbar{\Omega}$ and finally its normal spectrum is equal to $\partial\Omega.$ Consequently, the operator $M$ on $\big( H^2(\Omega), \lambda ds\big)$ must be unitarily equivalent to the bundle shift $T_{\alpha}$ on $\big( H^2_{\alpha}(\Omega), ds\big)$ for some character $\alpha.$ Now, we  compute the character $\alpha.$

Since $\lambda$ is a positive continuous function on $\partial \Omega,$ using Lemma $\ref{existence of modulus auto},$ we have the existence of a character $\alpha$  and a function $F$ in $H^{\infty}_{\alpha}(\Omega)$ satisfying $|F|^2 = \lambda$ on $\partial \Omega.$ The function $F$ is also invertible in the sense that there exist a function $G$ in  $H^{\infty}_{\alpha^{-1}}(\Omega)$ such that $FG=1$ on $\Omega.$ It is straightforward to verify that the linear map $M_F:\big(H^2(\Omega),\lambda ds\big)\mapsto \big(H^2_{\alpha}(\Omega),ds\big)$ defined by 
\begin{align*}
 M_F(g) &=F g, \;\;\;\;\;g\in \big(H^2 (\Omega),\lambda ds\big)
\end{align*}
is unitary. Also $M_F$ being a multiplication operator,  intertwines the corresponding multiplication operator by the coordinate function on the respective Hilbert spaces. From Lemma $\ref{existence of modulus auto},$ it is  clear that the character $\alpha$ is determined by the following $n$ tuple of numbers:
\begin{align}\label{cj(lambda)}
   c_j(\lambda) = -\int _{\partial \Omega _{j}} \frac{\partial}{\partial \eta _z}\big(u_{\lambda}(z)\big)ds(z), \;\;\;\mbox{for}\; j= 1,2,...,n,\;\;
\end{align} where $u_{\lambda}$ is the harmonic function on $\Omega$ with continuous boundary value $\frac{1}{2}\log \lambda.$ Using this information along with the Theorem $\ref{Equivalence of character}$, we deduce the following Lemma which describe the unitary equivalence class of the multiplication operator $M$ on $\big(H^2(\Omega),\lambda ds\big).$   
\begin{lem}\label{equivalence of measure}
Let $\lambda, \mu$ be two positive continuous function on $\partial \Omega.$ Then the operators $M$ on the Hilbert spaces $\big(H^2(\Omega),\lambda ds\big)$ and $\big(H^2(\Omega),\mu ds\big)$ are unitarily equivalent iff 
\begin{align*}
\exp\big(ic_{j}(\lambda)\big)&=\exp\big(i c_j(\mu)\big),\;\;\;  j= 1,\ldots, n.
\end{align*}
\end{lem}
It also follows from a result of Abrahamse (See \cite[Proposition 1.15]{ABtoeplitz}) that given a character $\alpha$ there exist a invertible element $F$ in $H^{\infty}_{\alpha}(\Omega)$ such that 
\begin{align*}
|F(z)|^2 &= \begin{cases}
    1,& \text{if } z\in \partial \Omega_{n+1}\\
    p_j,& \text{if } z\in \partial \Omega_{j},\;\;j=1,\cdots,n,             
\end{cases}
\end{align*} where $p_j$ are positive constant. Thus we have proved the following theorem.
\begin{thm}
Given any character $\alpha,$ there exists a positive continuous function $\lambda$ defined on $\partial\Omega$ such that the operator $M$ on $\big(H^2(\Omega),\lambda ds\big)$ is unitarily equivalent to the bundle shift $T_{\alpha}$ on $\big(H^2_{\alpha}(\Omega), ds\big).$
\end{thm}
 
\section{Weighted Kernel and Extremal Operator at a fixed point} 
Let $\lambda$ be a positive continuous function on $\partial\Omega.$ Since $\big( H^2(\Omega), d\omega_p\big)$ is a reproducing kernel Hilbert space and the norm on $\big( H^2(\Omega), d\omega_p\big)$ is equivalent to the norm on $\big( H^2(\Omega), \lambda ds\big),$ it follows that  $\big( H^2(\Omega), \lambda ds\big)$ is also a reproducing kernel Hilbert space. Let $K^{(\lambda)}(z,w)$ denote the kernel function for $\big( H^2(\Omega), \lambda ds\big).$

The case $\lambda \equiv 1$ gives us the S\"zego kernel $S(z,w)$ for the domain $\Omega.$ Associated to the S\"zego kernel, there exists a conjugate kernel $L(z,w),$ called the Garabedian kernel, which is related to the S\"{ze}ego kernel via the following identity. 
\begin{align*}
\overline{S(z,w)} ds &= \frac{1}{i} L(z,w) dz,\;\;w\in \Omega\;\text{and }\; z\in \partial \Omega
\end{align*}
We recall several well known properties of these two kernels when $\partial \Omega$  consists of jordan analytic curves. For each fixed $w$ in $\Omega$, the function $S_w(z)$ is holomorphic in a  neighbourhood of $\Omega$ and $L_w(z)$ is holomorphic in a neighbourhood of $\Omega -\{w\}$ with a simple pole at $w.$ $L_w(z)$ is non vanishing on $\overbar{\Omega} -\{w\}.$ The function $S_w(z)$ is non vanishing on $\partial \Omega$ and has exactly $n$ zero in $\Omega$ (See \cite[Theorem 13.1]{BELLcauchy}). In \cite[Theorem 1]{NEHARIexremal} Nehari has extended these result for the kernel $K^{(\lambda)}(z,w).$ 
\begin{thm}[Nehari]\label{Nehari}
Let $\Omega$ be a bounded domain in the complex plane, whose boundary consists of $n+1$ analytic jordan curve and let $\lambda$ be a positive continuous function on $\partial \Omega.$ Then there exist two analytic function $K^{(\lambda)}(z,w)$ and $L^{(\lambda)}(z,w)$ with the following properties: for each fixed $w$ in $\Omega$, the function $K^{(\lambda)}_w(z)$ and $L^{(\lambda)}_w(z) -(2\pi(z-w))^{-1}$ are holomorphic in $\Omega;$ $|K^{(\lambda)}_w(z)|$ is continuous on $\overbar{\Omega}$ and  $|L^{(\lambda)}_w(z)|$ is continuous in $\overbar{\Omega}- C_{\epsilon},$ where $C_{\epsilon}$ denotes a small open disc about $w$; $K^{(\lambda)}_w(z)$ and $L^{(\lambda)}_w(z)$ are connected by the identity
\begin{align}\label{eq:conjugate id 2}
\overline{K^{(\lambda)}_{w}(z)} \lambda(z) ds &= \frac{1}{i} L^{(\lambda)}_{w}(z) dz,\;\;w\in \Omega\;\text{and}\; z\in \partial \Omega 
\end{align}
These properties determine both functions uniquely. 
\end{thm}
From $(\ref{eq:conjugate id 2}),$ we have that $\frac{1}{i}K^{\lambda}_w(z) L^{\lambda}_w(z)dz \geq 0.$ The boundary $\partial\Omega$ consists of Jordan analytic curves, therefore from the Schwartz reflection principle, it follows that the function  $K^{\lambda}_w$ and $L^{\lambda}_w -(2\pi(z-w))^{-1}$ are holomorphic in a neighbourhood of $\overbar{\Omega}.$ 

We have shown that the operator $M$ on $\big( H^2(\Omega), \lambda ds\big)$ is unitarily equivalent to a bundle shift of rank $1$. Consequently the adjoint operator $M^*$ lies in $B_1(\Omega ^*)$ admitting $\overbar{\Omega^*}$ as a spectral set from which a curvature inequality follows: 

\begin{align*}
\mathcal{K}_{T} (w) &\leq -4\pi ^2(S_{\Omega ^*}(w,w))^2,\;w\, \in \Omega ^*.
\end{align*}
Or equivalently,
\begin{align*}
\frac{\partial ^2}{\partial w \partial \bar{w}}\mbox{log} K^{(\lambda)}(w,w)\geq 4\pi ^2 (S_{\Omega }(w,w))^2, \;\; \;\;w \in \Omega.
\end{align*}
Fix a  point $\zeta$ in $\Omega.$ The following lemma provides a criterion  for the adjoint operator $M^*$ on $\big( H^2(\Omega), \lambda ds\big)$ to be extremal  at $\bar{\zeta},$ that is,  
\begin{align*}
\frac{\partial ^2}{\partial w \partial \bar{w}}\mbox{log} K^{(\lambda)}(w,w) \mid _{w = \zeta} &= 4\pi ^2 (S_{\Omega }(\zeta,\zeta))^2.
\end{align*}
\begin{lem} \label{extremecriterion}
The operator $M^*$ on the Hilbert space $\big(H^2(\Omega), \lambda ds\big)$ is extremal at $\bar{\zeta}$ iff $L^{(\lambda)}_{\zeta}(z)$  and the S\"{z}ego kernel at $\zeta$, namely $S_{\zeta}(z)$ have the same set of zeros in $\Omega.$ 
\end{lem}
\begin{proof}
Consider the closed convex set  $M_1$ in $\big(H^2(\Omega), \lambda(z) ds\big)$ defined by
\begin{align*}
M_1 &:= \{f\in \big(H^2(\Omega), \lambda(z) ds\big): f(\zeta) =0, f'(\zeta) =1\}. 
\end{align*}
Now consider the extremal problem is 
\begin{align}\label{extremalprob}
\inf\; \{||f||^2:f\in M_1\}
\end{align}
Since $M_1$ is a closed convex set, there exist a unique function $F$ in $M_1$ which solve the extremal problem. It has been shown in \cite{GMCI} that the function $F$ in $\big(H^2(\Omega), \lambda(z) ds\big)$ is a solution to the extremal problem iff $F \in M_1$ and $F$ is orthogonal to the subspace
\begin{align*}
H_1 &= \{f\in \big(H^2(\Omega), \lambda(z) ds\big): f(\zeta) =0, f'(\zeta) =0 \} = \big(\rm Span \{K^{(\lambda)}_{\zeta}, \bar{\partial}K^{(\lambda)}_{\zeta}\}\big)^{\perp}.
\end{align*} 
Solution to this extremal problem can be found in terms of the Kernel function as in \cite{GMCI}:
\begin{align*}
\inf\; \{||f||^2:f\in M_1\} &= \bigg\{K^{(\lambda)}(\zeta,\zeta) \bigg(\frac{\partial ^2}{\partial w \partial \bar{w}}\mbox{log} K^{(\lambda)}(w,w)|_{w=\zeta}\bigg)  \bigg\}^{-1}.
\end{align*}
Now consider the function $g$ in $\big(H^2(\Omega), \lambda(z) ds\big)$ defined by
\begin{align*}
g(z):= \frac{ K^{(\lambda)}_{\zeta}(z) F_{\zeta}(z)}{2 \pi S(\zeta,\zeta) K^{(\lambda)}(\zeta,\zeta)}, \;\;\;z\in \Omega,
\end{align*} 
where $F_{\zeta}(z) = \frac{S_{\zeta}(z)}{L_{\zeta}(z)}$ denote the Ahlfors map for the domain $\Omega$ at the point $\zeta$ (See \cite[Theorem 13.1]{BELLcauchy}). Using the reproducing property for the kernel function  $K^{(\lambda)}$ and the fact that $|F_{\zeta} (z)| \equiv 1$ on $\partial \Omega$, it is straightforward to verify that 
\begin{align*}
 ||g||_{\lambda ds}^2 = \bigg( K^{(\lambda)}(\zeta,\zeta) 4\pi ^2 S(\zeta,\zeta)^2\bigg)^{-1}.
\end{align*}
Since $F_{\zeta}(\zeta)=0$ and $F_{\zeta}^\prime(\zeta)= 2 \pi S(\zeta,\zeta)$, it follows that $g\in M_1.$ Consequently we have 
\begin{align*}
\bigg( K^{(\lambda)}(\zeta,\zeta) 4\pi ^2 S(\zeta,\zeta)^2\bigg)^{-1} \geq \bigg\{K^{(\lambda)}(\zeta,\zeta) \bigg(\frac{\partial ^2}{\partial w \partial \bar{w}}\mbox{log} K^{(\lambda)}(w,w)|_{w=\zeta}\bigg)  \bigg\}^{-1}
\end{align*} Or equivalently
\begin{align*}
\frac{\partial ^2}{\partial w \partial \bar{w}}\mbox{log} K^{(\lambda)}(w,w) \mid _{w = \zeta} \geq 4\pi ^2 (S_{\Omega }(\zeta,\zeta))^2.
\end{align*}
So equality holds iff $g$ solve the extremal problem in $(\ref{extremalprob})$ iff  $g$ is orthogonal to the subspace $H_1.$ Hence, we conclude that the operator $M^*$ on the Hilbert space $\big(H^2(\Omega), \lambda(z) ds\big)$ is extremal at $\bar{\zeta}$  iff $g$ is orthogonal to the subspace $H_1$. Now consider the following integral
\begin{align*}
I_f &= \int _{\partial \Omega} f(z) \overline{K^{(\lambda)}_{\zeta}(z)} \overline{F_{\zeta}(z)} \lambda(z) ds\\
&= \frac{1}{i}\int _{\partial \Omega} f(z)\overline{F_{\zeta}(z)} L^{(\lambda)}_{\zeta}(z) dz\;\;\;\;\;\;\;\;(\text{Using the identity}\;\ref{eq:conjugate id 2})\\
&= \frac{2\pi}{2\pi i}\int _{\partial \Omega} \frac{f(z)}{F_{\zeta}(z)}L^{(\lambda)}_{\zeta}(z)dz\\
&= \frac{1}{2\pi i}\int _{\partial \Omega} \frac{f(z)L^{(\lambda)}_{\zeta}(z)\big(2\pi L_{\zeta}(z)\big)}{S_{\zeta}(z)}dz 
\end{align*}
Since $H_1 \cap \rm Rat(\overbar{\Omega})$ is dense in $H_1$,  $g$ is orthogonal to $H_1$ iff $I_f$ vanishes for all $f\in H_1 \cap \rm Rat(\overbar{\Omega})$. Observe that we have  $L^{(\lambda)}_{\zeta}(z) L_{\zeta}(z)$ is holomorphic in $\Omega-\{\zeta\}$ with a pole of order $2$ at $\zeta$. As $\partial \Omega$ consists of Jordan analytic curve,  both the function $L^{(\lambda)}_(\zeta)(z)$ and $L_{\zeta}(z)$ are also holomorphic in a neighbourhood of $\partial\Omega.$ It is known that $L_{\zeta}(z)$ has no zero in $\overbar{\Omega}-\{\zeta\}$ and $S_{\zeta}(z)$ has exactly $n$ zero say $a_1,a_2,...,a_n$ in $\Omega,$ (cf. \cite[]{BELLcauchy}).

Now we claim that $I_f$ vanishes for all $f\in H_1 \cap \rm Rat(\overbar{\Omega})$ iff the set of zeros of the function  $L^{(\lambda)}_{\zeta}(z)$ in $\Omega$ is $\{a_1,a_2,\ldots,a_n\}.$  

First if we assume that $L^{(\lambda)}_{\zeta}(z)$ has $\{a_1,a_2,\ldots,a_n\}$ as the zero set in $\Omega$, then the integrand in $I_f$ is holomorphic in a neighbourhood of $\overbar{\Omega}$ for every $f$ in $H_1 \cap \rm Rat(\overbar{\Omega})$ and consequently $I_f$ vanishes for every $f$ in $H_1 \cap \rm Rat(\overbar{\Omega}).$ Conversely if $L^{(\lambda)}_{\zeta}(z)$ doesn't vanish at one of $a_j$'s, without loss of generality, say at $a_1$, then the function
 \begin{align*}
 f =(z-\zeta)^2 \prod\limits_{k=2}^n (z-a_k) 
 \end{align*}
is in $H_1 \cap \rm Rat(\overbar{\Omega}).$ 
Observe that the integrand in $I_f,$ with  this choice of the function $f,$ is holomorphic in a neighbourhood of $\overbar{\Omega}$  except at the point $a_1,$ where it has a simple pole. So by the Residue theorem the integral $I_f$ equals  the residue of the integrand at $a_1,$ which is not zero completing the proof.  
\end{proof}

\subsection{Existence of Extremal operator}
We provide below two different descriptions of an  extremal operator at $\bar{\zeta}$ using the criterion obtained in Lemma \ref{extremecriterion}
Let $a_1,a_2,...,a_n$ be the zeros of the S\"{z}ego $S_{\zeta}(z)$ in $\Omega.$ 
\subsubsection{Realization of the extremal operator at $\bar{\zeta}$}
Consider the function $\lambda$ on $\partial \Omega$ defined by
\begin{align*}
\lambda(z) &:= \prod\limits_{k=1}^n |z-a_k|^2,\;\;z\in \partial\Omega.
\end{align*}
Then, for $z\in \partial \Omega,$ we have 
\begin{align*}
\frac{\overline{S_{\zeta}(z)}}{\prod \limits _{j=1}^n (\bar{z}-\bar{a}_j)(\zeta - a_j)} \lambda(z) ds &= \frac{\prod\limits_{k=1}^n (z-a_k)}{\prod\limits_{k=1}^n (\zeta -a_k)} \overline{S_{\zeta}(z)}ds\\
%\frac{\overline{S_{\zeta}(z)}}{\prod \limits _{j=1}^n (\bar{z}-\bar{a}_j)(\zeta - a_j)} \lambda(z) ds 
&= \frac{1}{i}\frac{\prod\limits_{k=1}^n (z-a_k)}{\prod\limits_{k=1}^n (\zeta -a_k)}L_{\zeta}(z) dz
\end{align*}
Note that the function $S_{\zeta}(z)\bigg( \prod \limits _{j=1}^n (z-a_j)(\bar{\zeta} - \bar{a}_j)\bigg)^{-1}$ is holomorphic in a neighborhood of $\overbar{\Omega}$ and the function $L_{\zeta}(z)\bigg( \prod\limits_{k=1}^n (z-a_k)\bigg) \bigg( \prod\limits_{k=1}^n (\zeta-a_k)\bigg)^{-1}$ is a meromorphic in a neighbourhood of $\overbar{\Omega}$ with a simple pole at $\zeta$. Hence using the uniqueness portion of the Theorem ($\ref{Nehari}$), we get
\begin{align*}
K^{(\lambda)}_{\zeta}(z) &= \frac{S_{\zeta}(z)}{\prod \limits _{j=1}^n (z-a_j)(\bar{\zeta} - \bar{a}_j)},\;z\in \overbar{\Omega} \;\;\text{and}\;\; L^{(\lambda)}_{\zeta}(z) = \frac{\prod\limits_{k=1}^n (z-a_k)}{\prod\limits_{k=1}^n (\zeta -a_k)}L_{\zeta}(z),\;\;z\in \overbar{\Omega}-\{\zeta\}.
\end{align*}
Clearly, $\{a_1,a_2,...,a_n\}$ is the zero set of the function $L^{(\lambda)}_{\zeta}(z).$ So, the adjoint operator $M^*$ on  $\big(H^2(\Omega), \lambda(z) ds\big)$ is an extremal operator  at $\bar{\zeta}.$

\subsubsection{A second realization of the extremal operator at $\bar{\zeta}:$} This realization of the extremal operator  was obtained earlier in \cite{GMCI}. Consider the measure 
\begin{align*}
\lambda(z) ds &= \frac{|S_{\zeta}(z)|^2}{S(\zeta,\zeta)}ds,\,\, z \in \partial \Omega,
\end{align*}
on the boundary $\partial \Omega.$ Using the reproducing property of the S\"zego kernel, it is easy to verify that 
\begin{align*}
\big\langle f, 1\big \rangle _{\big(H^2(\Omega), \lambda ds\big)} &= f(\zeta)
\end{align*}
This gives us $K^{(\lambda)}_{\zeta}(z)= 1$ for all $z\in\overbar{\Omega}.$ So we have
\begin{align*}
\lambda(z) ds &= \frac{S_{\zeta}(z)}{S(\zeta,\zeta)} \overline{S_{\zeta}(z)} ds,\;\;z\in\partial\Omega\\
&= \frac{1}{i}\frac{S_{\zeta}(z)}{S(\zeta,\zeta)} L_{\zeta}(z) dz,\;\;\;z\in\partial\Omega
\end{align*}
Now the function $S_{\zeta}(z)L_{\zeta}(z)\big( S(\zeta,\zeta)\big)^{-1}$ is a meromorphic function in a neighbourhood of $\overline{\Omega}$ with a simple pole at $\zeta$. Again, using the uniqueness guaranteed in  Theorem ($\ref{Nehari}$), we get
\begin{align*}
L^{(\lambda)}_{\zeta}(z) &= S_{\zeta}(z)L_{\zeta}(z)\big( S(\zeta,\zeta)\big)^{-1},\;\;\;\;\;\;z\in \overbar{\Omega}-\{\zeta\}.
\end{align*}
Again, the zero set of the function  $L^{(\lambda)}_{\zeta}(z)$ is $\{a_1,a_2,...,a_n\}.$ So the operator $M^*$ on $\big(H^2(\Omega), \lambda(z) ds\big)$ is an extremal operator at $\bar{\zeta}$.

We shall prove that the any extremal operator which is also the adjoint of a bundle shift is uniquely determined up to unitary equivalence. An amusing consequence of this uniqueness is that the two realizations of the extremal operators given above must coincide up to unitary equivalence. 

\subsection{Index of the Blaschke product} 
 
To facilitate the proof of the uniqueness, we need to recall basic properties of multiplicative Blaschke product on $\Omega$ and its index of automorphy. This is also going to be a crucial ingredient in  determining  the character $\alpha$ of the extremal operator at $\bar{\zeta}.$
 
%Recall the definition of Multiplicative Blaschke function on the domain $\Omega$. 

Let $g(z,a)$ be the Green's function for the domain $\Omega,$ whose critical point is $a \in \Omega$. The multiplicative Blaschke factor with zero at $a$, is defined as follows
   \begin{align*}
   B_{a}(z) &= \mbox{exp}(-g(z,a)-i g^*(z,a)),\;\;\mbox{for all}\;z\in \Omega,
   \end{align*} where $g^*(z,a)$ is the multivalued conjugate of the Green's function $g(z,a),$ which is harmonic on $\Omega - \{a\}.$ So, $B_a(z)$ is a multiplicative function on $\Omega,$ which vanishes only at the point $a$ with multiplicity $1$ and on $\partial \Omega$ its absolute value is identically $1$. Note that periods of the conjugate harmonic function $g^*(z,a)$ around the boundary component $\partial \Omega_j$ is equal to 
   \begin{eqnarray*}
   p_j(a) = -\int _{\partial \Omega _{j}} \frac{\partial}{\partial \eta _z}\big(g(z,a)\big)ds_z, \;\;\;\mbox{for}\; j= 1,2,...,n\;\;
   \end{eqnarray*}
The negative sign appearing in the equation for the periods is a result of the  assumption that $\partial \Omega$ is positively oriented, that is,  the boundary $\partial \Omega_j,$ $j=1,2,\ldots, n,$ except the outer one are oriented in clockwise direction.  
   
% Remember, If $\Gamma_{j}$ is positively oriented that is oriented in anti-clockwise direction, then $\frac{\partial}{\partial \eta _z} $ denote the directional derivative along the outward normal direction (w.r.t oriented $\Gamma_{j}$). Here we have assumed that $\Gamma$, the boundary of $\Omega$, is positively oriented. So for $j=1,2,...,n$, $\Gamma_{j}$ is oriented in clockwise direction and the outer boundary component $\Gamma_{n+1}$ is oriented in anticlockwise direction. Hence in the above expression and throughout the section, $\frac{\partial}{\partial \eta _z} $ denote the directional derivative along the outward normal direction (w.r.t oriented $\Gamma$). 
   
Since the Blaschke factor $B_a(z)$ is multiplicative function on $\Omega,$ therefore it is induced by a modulus automorphic function on unit disc, say $b_\alpha,$ for some $\alpha.$ The character $\alpha$ uniquely determines $n$-tuple of complex number of unit modulus. These are the image under $\alpha$ of the generators of the group $G$ of Deck transformations relative to the covering map $\pi:\mathbb D \to \Omega.$ This $n$ - tuple,  called the index of the Blaschke factor $B_a(z),$  is of the form   
 \begin{align*}
 \{\mbox{exp}(-i p_1(a)),\mbox{exp}(-i p_2(a),\cdots, \mbox{exp}(-i p_n(a))\} 
 \end{align*}
We recall bellow the well known relationship of the period $p_j(a)$  to the harmonic measure $\omega _j(z)$ of the boundary component $\partial\Omega_{j},$ namely, 
 \begin{align*}
 \omega _j(a) = -\frac{1}{2 \pi}\int  _{\partial \Omega _{j}} \frac{\partial}{\partial \eta _z}\big(g(z,a)\big)ds_z =\frac{1}{2 \pi}p_j(a) , \;\;\;\mbox{for}\; j= 1,2,...,n,
\end{align*}  
 where the harmonic measure $\omega _j(z)$ is the function which is harmonic in $\Omega$ and has the boundary values $1$ on $\partial\Omega_j$  and is $0$ on all the other boundary components. Hence the index of the Blaschke factor $B_a(z)$ is 
 \begin{align*}
 \mbox{ind}(B_a(z)) &=\{\mbox{exp}(-2\pi i \omega_1(a)),\mbox{exp}(-2 \pi i \omega_2(a),\cdots, \mbox{exp}(-2\pi i \omega_n(a))\}
 \end{align*}
 For each of these $n$ tuple of numbers, there exist a homomorphism $\alpha :G \to \mathbb{T}$ such that these $n$ tuple of numbers occur as  the image of the $n$ generator of the group $G$ under the map $\alpha$ completing the bijective correspondence between the character $\alpha$ and the index.  It follows that the function $B_a:=b_\alpha \circ \pi^{-1}$ lies in $H^{\infty}_{\alpha}.$ 
 
The index of the Blaschke product $B(z)=\prod \limits _{k=1}^m B_{a_k}(z),$  $a_k \in \Omega,$ is equal to 
 \begin{align}\label{eq:index of blaschke}
 \mbox{ind}(B(z)) &=\bigg\{\mbox{exp}\big(-2\pi i \sum\limits_{k=1}^{m}\omega_1(a_k)\big),\cdots, \mbox{exp}\big(-2\pi i\sum\limits_{k=1}^{m} \omega_n(a_k)\big)\bigg\}
\end{align}  

\subsection{Zeros of the S\"{z}ego kernel $S_\zeta(z)$} Fixing $\zeta$ in $\Omega,$ which is $n+1$ - connected, as pointed out earlier, the S\"{z}ego kernel $S_{\zeta}(z)$ has exactly $n$ zeros (counting multiplicity) in $\Omega.$ Let $a_1,a_2,\ldots, a_n$ be the zeros of $S_{\zeta}(z)$. Hence the Ahlfors function $F_{\zeta}(z)$ at the point $\zeta$ has exactly $n+1$ zeros in $\Omega,$ namely $\zeta,a_1,a_2,\ldots, a_n.$ Now an interesting relation between the  points $a_1, \ldots ,a_n$ and $\zeta$ becomes evident.

First consider the Blaschke product $B(z)= B_{\zeta}(z).\prod \limits _{k=1}^n B_{a_k}(z).$  The index of the Blaschke product $B(z),$ using $\eqref{eq:index of blaschke},$ is easily seen to be of the form
\begin{align*}
\beta = (\beta _1,\beta_2,\cdots,\beta_n), \mbox{ where } \beta_j &=\bigg\{\mbox{exp}\bigg(-2\pi i \big(\omega_j(\zeta)+\sum\limits_{k=1}^{n}\omega_j(a_k)\big )\bigg) \bigg\}\;\;\mbox{for}\;j=1,2,\ldots ,n.
\end{align*} 

The Ahlfors function  $F_{\zeta}(z)$ is in $H^\infty(\Omega)$ and it is holomorphic in a neighbourhood of $\bar{\Omega}$ as long as the boundary $\partial \Omega$ is analytic. 
Therefore in the inner outer factorization of $F_{\zeta}(z),$ there is no singular inner function and it follows that 
\begin{align*}
|F_{\zeta}(z)| &= |B(z)||\psi (z)|,\;\;z\in \Omega, 
\end{align*}
 where $\psi(z)$ is a multiplicative outer function of index 
 $$\beta ^{-1} = (\beta _1^{-1},\beta_2^{-1},\cdots,\beta_n^{-1}).$$ 
 
Now consider the linear map $ L:\big(H^2(\Omega),ds(z)\big)\mapsto \big(H^2_{\beta ^{-1} }(\Omega),ds\big),$ defined by
 \begin{align*}
Lf =\psi f, \;\;\;\;f\in \big(H^2(\Omega),ds(z)\big).
\end{align*}
Note that $\psi(z)$ is outer and it is bounded in absolute value (since $F_{\zeta}(z)$ is bounded) on $\Omega.$  It is straightforward to verify that $L$ is a unitary operator. Also, since  $L$ is a multiplication operator, it intertwines any two  multiplication operators on the respective Hilbert spaces. 

As a corollary of  Theorem \ref{Equivalence of character}, we must have $\beta ^{-1} = (1,1,\cdots,1).$ This implies 
\begin{align}\label{eq:invariance in zeros of Szego}
\mbox{exp}\bigg(-2\pi i \big(\omega_j(\zeta)+\sum\limits_{k=1}^{n}\omega_j(a_k)\big )\bigg) &= 1,\;\;\;j=1,2,\ldots,n,
\end{align}   
relating the point $\zeta$ to the zeros  $a_1, a_2,\cdots,a_n$ of the S\"{z}ego kernel $S_\zeta(z).$
\subsection{Uniqueness of the Extremal operator}

%We  have shown that for a positive continuous function $\lambda$ on $\partial \Omega,$  the operator $M^*$ on the Hilbert space $\big(H^2(\Omega), \lambda(z) ds\big)$ is extremal at $\bar{\zeta}$ iff $L^{(\lambda)}_{\zeta}(z)$ has exactly $\{a_1,a_2,...,a_n\}$ as the zero set; where $\{a_1,a_2,...,a_n\}$ are the zeros of the $S_{\zeta}(z)$ in $\Omega.$
 
Assume that for a positive continuous function $\lambda$ on $\partial\Omega,$ the operator $M^*$ on the Hilbert space $\big(H^2(\Omega), \lambda(z) ds\big)$ is extremal at $\bar{\zeta}$. The function $K^{(\lambda)}_{\zeta}(z)$ is analytic in a neighborhood of $\overline{\Omega}$ and the conjugate kernel $L^{(\lambda)}_{\zeta}(z)$ is meromorphic in a neighborhood of $\overline{\Omega}$ with a simple pole only at the point $\zeta$. Also from Lemma \ref{extremecriterion}, we have that the zero set of $L^{(\lambda)}_{\zeta}(z)$ is the the set  $\{a_1,a_2,...,a_n\},$ where $\{a_1,a_2,...,a_n\}$ are the zeros of  $S_{\zeta}(z)$ in $\Omega.$ We have, using the equation $\eqref{eq:conjugate id 2}$, that
 \begin{align*}
|K^{(\lambda)}_{\zeta}(z)|^2 \lambda(z)ds &= \frac{1}{i}K^{(\lambda)}_{\zeta}(z)L^{(\lambda)}_{\zeta}(z)dz,\;\;\;\;z\in \partial\Omega.
\end{align*}
An application of the Generalized  Argument Principle shows that the total number of zeros of $K^{(\lambda)}_{\zeta}(z)$ and $L^{(\lambda)}_{\zeta}(z)$  in $\overbar{\Omega}$, where a zero on the boundary is counted as $\tfrac{1}{2}$, is equal to $n$. Hence it follows that  $a_1,a_2,...,a_n$ are the all zeros of $L^{(\lambda)}_{\zeta}(z)$ in $\overbar{\Omega}$ and $K^{(\lambda)}_{\zeta}(z)$ has no zero in $\overbar{\Omega}.$

Nehari  \cite[Theorem 4]{NEHARIexremal}  has shown that the meromorphic function 
\begin{align*}
R(z) &= \frac{K^{(\lambda)}_{\zeta}(z)}{L^{(\lambda)}_{\zeta}(z)},\;\;\;\;z\in \overbar{\Omega}
\end{align*}
with exactly one zero at $\zeta$ and poles exactly at $a_1,a_2,...,a_n$, solves the extremal problem 
\begin{align*}
\sup\{|f'(\zeta)| : f \in B_{\lambda}\},
\end{align*}
where $B_{\lambda}$ denotes the class of meromorphic function on $\Omega.$ Each $f$ in $B_\lambda$ is required to vanish at $\zeta$ and it is assumed that the set of poles of $f$ is a subset of $\{a_1,a_2,...,a_n\}.$  The radial limit of the functions $f$ at $z_0 \in \partial \Omega,$ from within $\Omega,$ in the class $B_\lambda$ are uniformly bounded:
\begin{align*}
\limsup_{z \to z_0} |f(z)| \leq \frac{1}{\lambda(z_0)},\;\;\;\;z_0\in \partial \Omega.
\end{align*} 
The proof includes the verification 
\begin{align*}
|R(z)| &= \frac{1}{\lambda(z)},\;\;z\in \partial \Omega
\end{align*}

Now consider the multiplicative function $G$ on $\Omega$ defined by 
\begin{align*}
G(z) &= \frac{B_{\zeta}(z)}{R(z)\prod \limits_{j=1}^n B_{a_j}(z)},\;\;\;\;\;\;z\in \overbar{\Omega}.
\end{align*}
So, $G$ is a multiplicative function in a neighbourhood of $\overbar{\Omega}.$ Also by construction $|G|$ has no zero in $\overbar{\Omega}.$ Using the inner outer factorization for multiplicative functions (See \cite[Theorem 1]{VOICHICKinner}), we see   that $G$ is a bounded multiplicative outer function. Also note that 
\begin{align*}
|G(z)|&=\lambda(z),\;\;\; z \in \partial \Omega.
\end{align*}
The index of $G$ is given by
\begin{align*}
\bigg\{\exp\bigg(2\pi i\big(-\omega _1(\zeta) +\sum\limits _{j=1}^n \omega _1(a_j)\big)\bigg),...,\exp\bigg(2\pi i\big(-\omega _n(\zeta) + \sum\limits _{j=1}^n \omega _n(a_j)\big)\bigg)\bigg\}
\end{align*}
Using equation $\eqref{eq:invariance in zeros of Szego}$, we infer that the index of $G(z)$ must be equal to
  \begin{align*}
\bigg\{\exp\bigg(-4\pi i\omega _1(\zeta)\bigg) ,...,\exp\bigg(-4\pi i\omega _n(\zeta)\bigg) \bigg\}.
\end{align*}

The function  $G$ is outer and hence the function $F:=\sqrt{G}$ is well defined.  It is a bounded multiplicative outer function with $|F(z)|^2 = \lambda(z)$ for all $z$ in $\partial\Omega.$ Let's denote the index of $F$  by 
\begin{align*}
\bigg\{\exp\bigg(-2\pi i\omega _1(\zeta)\bigg) ,...,\exp\bigg(-2\pi i\omega _n(\zeta)\bigg) \bigg\}.
\end{align*}

Now consider the linear map $V:\big(H^2(\Omega),\lambda(z)\;ds\big)\mapsto \big(H^2_{\alpha }(\Omega),ds\big)$ defined by
 \begin{align*}
Vf &=F f, \;\;\;\;\;f\in \big(H^2(\Omega),\lambda(z)\;ds\big).
\end{align*}
It is easily verified that $V$ is a unitary multiplication operator, which intertwines the corresponding multiplication operators on the respective Hilbert spaces. Hence  the character $\alpha$ for the bundle shift $T_{\alpha}$ on $\big(H^2_{\alpha}(\Omega),ds\big),$ which is extremal at $\bar{\zeta},$ is uniquely determined by the following $n$ tuple of complex number of unit modulus: 
\begin{align*}
\big\{\exp\big(-2\pi i\omega _1(\zeta)\big) ,...,\exp\big(-2\pi i\omega _n(\zeta)\big) \big\}
=\big\{\exp\big(2\pi i(1-\omega _1(\zeta))\big) ,...,\exp\big(2\pi i(1-\omega _n(\zeta))\big) \big\}.
\end{align*}
Hence if the adjoint of a bundle shift (upto unitary equivalence) is extremal at $\bar{\zeta},$ then it is uniquely determined. This completes the proof of the Theorem \ref{Uniqueness}. 

Since the group of the Deck transformations $G$ for the covering $\pi:\mathbb D \to \Omega$ is isomorphic to the free group on $n$ generators, 
any character $\alpha$ of the group $G$ is unambiguously determined, up to a permutation in the choice of generators for the group $G,$ by the $n$-tuple $\{x:= (x_1,x_2,\cdots,x_n): x_1, \ldots ,x_n \in [0,1)\},$ namely, 
$$ \alpha(g_k) = \exp( 2 \pi  i x_k), \,x_k\in [0,1)\, 1\leq k \leq n,$$
where $g_k,\,1\leq k \leq n,$ are generators of the group $G.$
The unitary equivalence class of the bundle shifts $T_\alpha$ of rank $1$ is therefore determined by the $n$-tuple $x$ in $[0,1)^n$ corresponding to the character $\alpha.$ 

%The set $\big \{ x:= (x_1,x_2,\ldots,x_n)\mid x_1, \ldots ,x_n \in [0,1)\big \}$ is homeomorphic to the $n$-torus, with the usual identification of the endpoints, via the map 
%$$x\mapsto \big(\exp(2\pi ix_1),\exp(2\pi i x_2),...,\exp(2 \pi ix_n)\big),$$
%and the class of the bundles shift $T_\alpha$ on $\big(H^2_{\alpha}(\Omega),ds\big)$ is uniquely determined by the character $\alpha,$ it follows that 

For $\zeta$ in $\Omega,$ the character corresponding to the $n$-tuple  $\bigg( (1-\omega _1(\zeta)), (1-\omega _2(\zeta)),\cdots,(1- \omega_n(\zeta))\bigg)$ defines the bundle shift which is extremal at $\bar{\zeta}.$ Let $\phi: \Omega \to [0,1)^n$ be the induced map, that is, 
\begin{align*}
\phi (\zeta) &= \bigg( (1-\omega _1(\zeta)), (1-\omega _2(\zeta)),\cdots,(1- \omega_n(\zeta))\bigg).
\end{align*}

Suita \cite{SUITA} shows that the map $\phi$ is not onto since $(0,\ldots ,0),$ which corresponds to the operator $M^*$ on the usual Hardy space, cannot be in its range. However, we show below that many other bundle shifts are missing from the range of the map $\phi,$ when $n\geq 2.$  

Let $\omega _{n+1}(z)$ be the harmonic measure for the outer boundary component $\partial \Omega_{n+1}.$ Thus $\omega _{n+1}$ is the harmonic function on $\Omega$ which is $1$ on $\partial\Omega_{n+1}$  and is $0$ on all the other boundary components. We have 
\begin{align*}
\sum\limits_{j=1}^{n+1}\omega _j \equiv 1 \;\;\;\mbox{and}\;\;\; 0< \omega_{n+1}(z) <1,\;\;z\in \Omega
\end{align*}
therefore 
\begin{align*}
(n-1) < \sum\limits_{j=1}^{n}\big(1-\omega _j(\zeta)\big) < n.
\end{align*}

From this, for $n\geq 2$, it follows that the set of extremal operators does not include the adjoint of many of the  bundle shifts. For instance, if the index of a bundle shift is $(x_1, \ldots, x_n)$ in $[0,1)^n$ is such that $x_1 + \cdots + x_n < n-1,$ then it cannot be an extremal operator at any $\bar{\zeta},$ $\zeta \in \Omega.$

\section{The special case of the Annulus}

Let $\Omega$ be an Annulur domain $A(0;R,1)$ with inner radius $R,$  $0<R<1,$ and outer radius $1.$ In this case we have a explicit expression for the harmonic measure corresponding to the boundary component $\partial \Omega _1,$ namely, 
\begin{align*}
\omega_1(z) = \frac{\log |z|}{\log R}.
\end{align*}  So for a fixed point $\zeta$ in $A(0;R,1),$ the character of the unique bundle shift which happens to be an extremal operator at $\bar{\zeta}$ is determined by the number 
 \begin{align*}
 \alpha(\zeta) = \exp \bigg(2\pi i\big(1- \omega_1(\zeta)\big)\bigg)
\end{align*}  
From this expression for the index, it is clear, in the case of an Annulur domain $A(0;R,1),$ that the adjoint of every bundle shift except the trivial one, is an extremal operator at some point $\bar{\zeta}$ in $\Omega ^*.$ In fact this is true of any doubly connected bounded domain $\Omega$ with Jordan analytic boundary since for such domain we have $\omega _1(\Omega)= (0,1),$ where $\omega_1$ is the harmonic measure corresponding to the inner boundary component $\partial \Omega _1.$ 

We now give a different proof  of the  Theorem \ref{Uniqueness} in the case of  $\Omega = A(0;R,1).$ In the course of this proof we see the effect of the weights on the zeros of the weighted Hardy kernels $K^{(\alpha)}.$ This question was raised in \cite{MCULAszego}.

For a fixed real number $\alpha$, Consider the measure $\mu_{\alpha}ds$ on the boundary of the Annulus, where the function $\mu_{\alpha}$ is defined by

\begin{equation*}
\mu_{\alpha}(z)=\begin{cases}
    1, & \text{if $|z|=1$},\\
    R^{2\alpha}, & \text{$|z|=R$}.
  \end{cases}
\end{equation*}
It is straightforward to verify that the function $\{f_n(z)\}_{n\in \mathbb{Z}}$ defined by 
\begin{align*}
f_n(z) &= \frac{z^n}{\sqrt{2\pi (1+R^{2\alpha +2n +1} )}},\;\;\;\;\;n\in \mathbb{Z}\
\end{align*}
forms an orthonormal basis for the Hilbert space $\big(H^2(\Omega), \mu_{\alpha} ds\big).$  The function 
\begin{align*}
K^{(\alpha)}(z,w) &:= \frac{1}{2 \pi} \sum\limits_{n\in \mathbb{Z}}\frac{(z\bar{w})^n}{1+R^{2\alpha +2n +1}},\;\;\;\;z,w \in \Omega,
\end{align*}
is uniformly convergent on compact subsets of $\Omega.$ Hence $K^{(\alpha)}$ 
is the reproducing kernel of the Hilbert space $\big(H^2(\Omega), \mu_{\alpha} ds\big).$ For each fixed $w$ in $\Omega,$ the kernel function $K^{(\alpha)}(z,w)$ is defined on $\Omega.$ However, it extends  analytically to a larger domain. To describe this extension,  recall that the Jordan Kronecker function, introduced by Venkatachaliengar (cf. \cite[p.37]{VENKATA}), is given by the formula 
\begin{align*}
f(b,t) &= \sum\limits_{k\in \mathbb{Z}}\frac{t^n}{1-bR^{2n}}. 
\end{align*}
This series converges for $R^2 < |t| <1,$ and for all $ b\neq R^{2k}, k\in\mathbb{Z}.$ Venkatachaliengar, using Ramanujan's $_{1}\psi_{1}$ summation formula,  has established the following identity (See \cite[p. 40]{VENKATA})
\begin{align}\label{Ramanujan id}
f(b,t) &=\frac{\prod\limits_{j=0}^{\infty}(1-btR^{2j})\prod\limits_{j=0}^{\infty}(1-\frac{R^{2j+2}}{bt})\prod\limits_{j=0}^{\infty}(1-R^{2j+2})\prod\limits_{j=0}^{\infty}(1-R^{2j+2})}{\prod\limits_{j=0}^{\infty}(1-tR^{2j})\prod\limits_{j=0}^{\infty}(1-\frac{R^{2j+2}}{t})\prod\limits_{j=0}^{\infty}(1-bR^{2j})\prod\limits_{j=0}^{\infty}(1-\frac{R^{2j+2}}{b})}.
\end{align}
This extends the definition of $f(b,t),$ as a meromorphic function, to all of the complex plane with a simple poles at $b=R^{2k},\, t=  R^{2k};\, k\in \mathbb{Z}.$ For fixed $w$ in $\Omega,$ since the function $f( -R^{2\alpha +1}, z\bar{w})$ coincides with $2\pi K^{(\alpha)}(z,w)$ for all $z$ in $\Omega,$ and $f$ is a meromorphic on the entire complex plane, it follows that $K_w^{(\alpha)}$ also extends to all of $\mathbb C$ as a meromorphic function. The poles of $K_w^{(\alpha)}$ are exactly at $\frac{R^{2k}}{\bar{w}}, k\in \mathbb{Z}.$ 
The zeros of the kernel function $K^{(\alpha)}_{w}(z)$ in $\Omega$ can  also be computed using the equation \ref{Ramanujan id}. The zeros $(b,t)$ of the function  $f$ must satisfy one of the following identities 
\begin{align*}
bt = R^{-2j}, j=0,1,2,\ldots \;\;\; \text{or,}\;\;\;\; bt = R^{2j+2}, j=0,1,2,\ldots 
\end{align*} 
For example, when $\alpha =0,$ the kernel $K^{(\alpha)}(z,w)$ is the S\"{z}ego kernel $S(z,w).$ It follows that if $w$ is a fixed but arbitrary point in $\Omega,$ then the zero set of the S\"{z}ego kernel function $S_w(z)$ is $\{-\frac{R}{\bar{w}}\}.$

The operator $M$ on the Hilbert space $\big(H^2(\Omega), \mu_{\alpha} ds\big)$ is a bilateral weighted shift with weight sequence 
\begin{align*}
\omega^{(\alpha)}_n &= \sqrt{\frac{1+ R^{2\alpha +2n +3}}{1+ R^{2\alpha +2n +1}}}, \;\;\;\;\;\;n\in \mathbb{Z}.  
\end{align*}
The identity 
\begin{align*}
\omega^{(\alpha +1)}_n &= \omega^{(\alpha)}_{n+1},\;\;n\in \mathbb{Z},  
\end{align*}
makes the operators $M$ on $\big(H^2(\Omega), \mu_{\alpha} ds\big)$ and $\big(H^2(\Omega), \mu_{\alpha +1} ds\big)$ unitarily equivalent.  Thus there is a natural map from the unitary equivalence classes of these bi-lateral shifts onto $[0,1).$  In the case of the annulus $A(0;R,1),$ we find that $u_{\mu_\alpha}$  and $c_{1}(\mu_\alpha),$ as defined in equation \eqref{cj(lambda)}, are equal to $\alpha \log |z|$ and $2\pi \alpha$ respectively. Applying Lemma $\ref{equivalence of measure}$, we see that  the operators $M$  on $\big(H^2(\Omega), \mu_{\alpha} ds\big)$ and $\big(H^2(\Omega), \mu_{\alpha +1} ds\big)$ are unitarily equivalent iff $\alpha - \beta$ is an integer. Thus we have a bijective correspondence between the unitary equivalence classes of these bi-lateral shifts and $[0,1),$ and we may assume  without loss of generality that $\alpha \in [0,1)$. 

For each $\alpha \in [0,1)$, the operator $M$  on $\big(H^2(\Omega), \mu_{\alpha} ds\big)$ is unitarily equivalent to the bundle shift $T_\beta$ on $\big(H^2_{\beta}(\Omega),  ds\big),$ where 
the character $\beta$ is determined  by the unimodular scalar $\exp(2\pi i\alpha).$

Now Fix a point $\zeta$ in $\Omega$. It is known that $S_{\zeta}(z)$, the S\"{z}ego kernel at $\zeta$ for the domain $\Omega$ has exactly one zero at $-\frac{R}{\bar{\zeta}}.$  
The existence of a conjugate kernel $L^{(\alpha)}(z,w)$ is established in \cite{NEHARIexremal}.  Then using the characterization for the extremal operator at $\bar{\zeta}$, it follows that the operator $M^*$  on $\big(H^2(\Omega), \mu_{\alpha} ds\big)$ is extremal at $\bar{\zeta}$ iff $ L^{(\alpha)}_{\zeta}(-\frac{R}{\bar{\zeta}})=0.$ 
From the identity 
\begin{align*}
zL^{(\alpha)}(z,w) &= K^{(\alpha)}(\tfrac{1}{z}, \bar{w})
\end{align*}
proved in \cite[p.1118]{MCULAszego}, and recalling that  $K^{(\alpha)}_{\bar{\zeta}}(-\frac{\bar{\zeta}}{R}) = \sum\limits_{k\in \mathbb{Z}}\frac{(-\frac{|\zeta|^2}{R})^n}{1+ R^{2\alpha +2n+1}},$ we conclude:  The operator $M^*$  is extremal at $\bar{\zeta}$ iff  
\begin{align*}
\sum\limits_{k\in \mathbb{Z}}\frac{(-\tfrac{|\zeta|^2}{R})^n}{1+ R^{2\alpha +2n+1}} = 0.
\end{align*}
Consequently, the operator $M^*$  on $\big(H^2(\Omega), \mu_{\alpha} ds\big)$ is extremal at $\bar{\zeta}$ iff the Jordan Kronecker function $f$ satisfy 
\begin{align*}
f(-R^{2\alpha +1},-\tfrac{|\zeta|^2}{R})&=0.
\end{align*}
So for a fixed $\zeta$, the real number $\alpha\in [0,1)$ must satisfy at least one of these identities 
\begin{align*}
R^{2\alpha} |\zeta|^2 = R^{-2j},j=0,1,2,...; \text{or}\;R^{2\alpha} |\zeta|^2 = R^{2j+2},j=0,1,2,\ldots 
\end{align*}
In any case, one must have 
\begin{align*}
\alpha &= \big(1-\frac{\log|\zeta|}{\log R}\big)\,\,(\text{mod } 1)
\end{align*}

So the unitary equivalence class of an operator which is extremal  at $\bar{\zeta}$, and is the adjoint of a bundle shift is uniquely determined.  Hence we have proved the Theorem stated below. 
\begin{thm}
The operator $M^*$  on $\big(H^2(\Omega), \mu_{\alpha} ds\big)$ is extremal at $\bar{\zeta}$ iff $\alpha = \big(1-\frac{\log|\zeta|}{\log R}\big)\,\,(\text{mod } 1).$
\end{thm}
%--------------------------------------------------------------------------
\subsection*{Acknowledgement} The author thanks G. Misra for his patient guidance and suggestions in the preparation of this paper. He would also like to thank the Math Stack Exchange community for providing an excellent opportunity for many stimulating discussions. 
\bibliographystyle{amsplain}\bibliography{bibfile}
\end{document}